\documentclass[a4paper,11pt,reqno]{amsart}
\usepackage{amsmath}
\usepackage{amsfonts}
\usepackage{amssymb}
\usepackage{amsthm}
\usepackage[all]{xy}
\usepackage[CJKbookmarks]{hyperref}
\usepackage{fancyhdr}
\usepackage{geometry}
\usepackage{color}

\theoremstyle{definition}
\newtheorem{definition}{Definition}[section]

\newtheorem{example}[definition]{Example}

\theoremstyle{plain}
\newtheorem{theorem}[definition]{Theorem}

\newtheorem{lemma}[definition]{Lemma}
\newtheorem{proposition}[definition]{Proposition}
\newtheorem{corollary}[definition]{Corollary}

\DeclareMathOperator{\cA}{\mathcal{A}}

\DeclareMathOperator{\cE}{\mathcal{E}}

\DeclareMathOperator{\cK}{\mathcal{K}}
\DeclareMathOperator{\cL}{\mathcal{L}}
\DeclareMathOperator{\cM}{\mathcal{M}}
\DeclareMathOperator{\cN}{\mathcal{N}}

\DeclareMathOperator{\cQ}{\mathcal{Q}}

\DeclareMathOperator{\cX}{\mathcal{X}}

\DeclareMathOperator{\bN}{\mathbb{N}}
\DeclareMathOperator{\bQ}{\mathbb{Q}}
\DeclareMathOperator{\bZ}{\mathbb{Z}}
\DeclareMathOperator{\fm}{\mathfrak{m}}
\DeclareMathOperator{\fp}{\mathfrak{p}}

\DeclareMathOperator{\Ann}{Ann}

\DeclareMathOperator{\bideg}{bideg}

\DeclareMathOperator{\Coker}{Coker}

\DeclareMathOperator{\coh}{coh}

\DeclareMathOperator{\Ext}{Ext}
\DeclareMathOperator{\gExt}{\underline{Ext}}

\DeclareMathOperator{\grKdim}{gr.Kdim}
\DeclareMathOperator{\grgldim}{gr.gldim}
\DeclareMathOperator{\grheight}{gr.ht}
\DeclareMathOperator{\grpdim}{gr.pdim}

\DeclareMathOperator{\gr}{gr}

\DeclareMathOperator{\gldim}{gldim}

\DeclareMathOperator{\Gr}{Gr}

\DeclareMathOperator{\grade}{grade}

\DeclareMathOperator{\Hom}{Hom}
\DeclareMathOperator{\gHom}{\underline{Hom}}
\DeclareMathOperator{\height}{ht}

\DeclareMathOperator{\idim}{inj.dim}

\DeclareMathOperator{\Kdim}{Kdim}

\DeclareMathOperator{\Lim}{Lim}

\DeclareMathOperator{\Mod}{Mod}

\DeclareMathOperator{\pdim}{pdim}

\DeclareMathOperator{\Proj}{Proj}
\DeclareMathOperator{\qgr}{qgr}
\DeclareMathOperator{\QGr}{QGr}
\DeclareMathOperator{\rank}{rank}

\DeclareMathOperator{\Spec}{Spec}

\DeclareMathOperator{\Tor}{Tor}
\DeclareMathOperator{\tor}{tor}

\newcommand{\dlim}{\underrightarrow{\Lim}}

\begin{document}
\title{Regular $\bZ$-graded local rings and  Graded Isolated Singularities}
\author{Haonan Li}

\author{Quanshui Wu}
\address{School of Mathematical Sciences \\
Fudan University \\
Shanghai, 200433 \\
 China}
\email{lihn@fudan.edu.cn, qswu@fudan.edu.cn}
\thanks{This research has been supported by the NSFC (Grant No. 12471032) and the National Key Research and Development Program of China (Grant No. 2020YFA0713200).}

%14A22 Noncommutative algebraic geometry [See also 16S38]
%16S38 Rings arising from noncommutative algebraic geometry [See also 14A22]
%16W50 Graded rings and modules (associative rings and algebras)
%16E35 Derived categories and associative algebras
%13D07 Homological functors on modules of commutative rings (Tor, Ext, etc.)
%13A02 Graded rings [See also 16W50]
%13H05 Regular local rings
%16E65 Homological conditions on associative rings (generalizations of regular, Gorenstein, Cohen-Macaulay rings, etc.)

\date{}

\newgeometry{left=3.18cm,right=3.18cm,top=2.54cm,bottom=2.54cm}

\begin{abstract}
  In this note we first study regular $\mathbb{Z}$-graded local rings. We characterize  commutative noetherian regular $\mathbb{Z}$-graded local rings in similar ways as in the usual local case.
    Then, we characterize graded isolated singularity for a commutative $\mathbb{Z}$-graded semilocal algebra in terms of the global dimension of its associated noncommutative projective scheme.
    As a corollary, we obtain that a commutative affine $\mathbb{N}$-graded algebra generated in degree $1$ is a graded isolated singularity if and only if its associated noncommutative projective scheme is smooth; if and only if the category of coherent sheaves on its projective scheme has finite global dimension, which are known in literature.
\end{abstract}

\maketitle

\section{Introduction}
A commutative noetherian local ring $(R,\fm)$ is called an isolated singularity if the scheme $\Spec R\backslash \{\fm\}$ is smooth, or equivalently, $R_{\fp}$ is a regular local ring for any non-maximal prime ideal $\fp$. For a commutative noetherian $\bN$-graded algebra $A$, $A$ is called a graded isolated singularity if the associated  projective scheme $\Proj A$
(an analogue of $\Spec R\backslash\{\fm\}$ in affine case) is smooth, that is, the degree-zero part of the homogeneous localization $A_{(\fp)}$ is regular for any graded prime ideal $\fp$ not containing $A_{\geqslant 1}$, which is equivalent to that the global dimension of $\coh (\Proj A)$  is finite, where $\coh (\Proj A)$ is the abelian category of the coherent sheaves on $\Proj A$.

Any commutative local Cohen-Macaulay ring of finite Cohen-Macaulay type is an isolated singularity (see \cite[Corollary 2]{HL},\cite[Theorem, p. 234]{Aus}). Motivated by \cite{HL}, a noncommutative analogue of isolated singularities for noncommutative connected graded algebras is considered in \cite{Jo}, and studied further in \cite{SV, Ue1, Ue2, MU} etc.

A well known result of Serre \cite{Se} says that for any commutative affine graded algebra $A$ generated in degree $0$ and $1$, the category $\coh (\Proj A)$ of coherent sheaves is equivalent to the quotient category $\qgr A$, the category of finitely generated graded $A$-modules modulo finite-dimensional $A$-modules. In fact, the quotient category $\qgr A$ is defined for any noetherian (not necessarily commutative) $\bN$-graded algebra $A$.
Inspired by Serre's result, $\qgr A$ is called the noncommutative projective scheme associated to $A$ in \cite{AZ}.  A noetherian $\bN$-graded algebra $A$ is called a noncommutative isolated singularity if $\qgr A$ has finite global dimension \cite{Ue1} (see Definition \ref{def-nc-iso-sing}).

There are some subtle places to be paid more attention in the proof of the above facts. For example, homogeneous localization $A_{(p)}$ is $\bZ$-graded, usually not bounded-below even $A$ is $\bN$-graded; Prime avoidance lemma does not hold in general in the graded case;  The relation between the Ext-groups $\Ext_{\qgr A}^i(\cM,\cN)$ in abelian category $\qgr A$ and  $\Ext_{\QGr A}^i(\cM,\cN)$ in Grothendieck category $\QGr A$ for any $\cM, \cN \in \qgr A$ should be clarified.

In this note, we first define regular $\bZ$-graded local rings. A noetherian commutative $\bZ$-graded local ring $(A,\fm)$ is called regular if its graded Krull dimension is equal to the rank of $\fm/{\fm^2}$ as $k_A$-module where $k_A=A/\fm$ (see Definition \ref{def-gr-reg-local-ring}). 
Regular graded local rings have similar characterizations as the usual regular local rings. 

\begin{theorem}[Theorem \ref{chara of graded regular}]\label{chara of graded regular in intro}
 Let $(A,\fm,k_A)$ be a noetherian $\bZ$-graded local ring of graded Krull dimension $d$. Then the following statements 
 are equivalent.
    \begin{itemize}
        \item [(1)] $(A,\fm)$ is a regular graded local ring.
        \item [(2)] $(A_{\fm},\fm A_{\fm})$ is a regular local ring.
       \item [(3)] The graded global dimension of $A $ is finite.
        \item [(4)] The graded global dimension of $A$ is equal to its graded Krull dimension.
        \item [(5)] $A$ is a regular ring.
        \item [(6)] $k_A[x_1,\cdots,x_d] \cong \Gr_{\fm} A$ as bigraded rings, where the bigrading of the polynomial ring $k_A[x_1,\cdots,x_d]$ is given as in the following: the first grading of $k_A$ is inherited from $k_A=A/\fm$, and the second grading of elements in $k_A$ are zero; the second grading of the homogeneous elements $x_i$ are all $1$.
        \item [(7)] $A$ has a regular sequence of the length $\rank_{k_A}(\fm/\fm^2)$ consisting of homogeneous elements.
    \end{itemize}
\end{theorem}

By using Theorem \ref{chara of graded regular in intro}, we give a detailed proof of the fact that a $\bZ$-graded (semi)local ring $A$ is a graded isolated singularity if and only if the global dimension of the abelian category $\qgr A$ is finite.
\begin{theorem}[Theorem \ref{char-graded-isolated-singualrty}]
    Let $A$ be a commutative noetherian $\bZ$-graded semilocal ring of graded Krull dimension $d$.
    Then the following are equivalent.
    \begin{itemize}
        \item [(1)] $(A_{(\fm)},\fm A_{(\fm)})$ is a graded isolated singularity for any maximal graded ideal $\fm$ of $A$.
        \item [(2)] $(A_{\fm},\fm A_{\fm})$ is an isolated singularity for any maximal graded ideal $\fm$ of $A$.
        \item [(3)] $(A_{(\fp)},\fp A_{(\fp)})$ is a regular graded local ring for any non-maximal graded prime ideal $\fp$ of $A$.
        \item [(4)] $(A_{\fp},\fp A_{\fp})$ is a regular local ring for any non-maximal graded prime ideal $\fp$ of $A$.
        \item [(5)] The global dimension of $\qgr A$ is $d-1$.
        \item [(6)] The global dimension of $\qgr A$ is finite.
    \end{itemize}
\end{theorem}

Then, for any graded quotient ring $A$ of a polynomial algebra, we obtain the characterization of $A$ being a graded isolated singularity in terms of $\Proj A$ and $\qgr A$ in literature, which also justifies the definition of noncommutative graded isolated singularities.

\begin{corollary}[Corollary \ref{iso-sing-of-quo-of-poly}]
    Let $A$ be a commutative affine graded algebra generated in degree $1$, $\fm=A_{>0}$. Let $\Proj A$ be the associated projective scheme of $A$. Then the following are equivalent.
\begin{itemize}
    \item [(1)] $(A,\fm)$ is a graded isolated singularity.
    \item [(2)] $(A_{\fm},\fm A_{\fm})$ is an isolated singularity.
    \item [(3)] The global dimension of $\qgr A$ is finite.
    \item [(4)] The global dimension of $\coh (\Proj A)$ is finite.
    \item [(5)] For any $\fp\in\Spec A\backslash\{\fm\}$, $(A_{\fp},\fp A_{\fp})$ is a regular local ring.
    \item [(6)] $\Proj A$ is smooth.
\end{itemize}
In this case, the global dimensions of $\qgr A$ and $\coh (\Proj A)$ are $\grKdim A-1$.
\end{corollary}

\section{Preliminaries}
\subsection*{Graded rings and modules}
A ring $A$ is called $\bZ$-graded if $A=\mathop{\oplus}\limits_{i\in \bZ} A_i$ where $\{A_i\mid i\in \bZ\}$ is a family of additive subgroups of $A$ such that $A_iA_j\subseteq A_{i+j}$ for all $i,j\in \bZ$. A $\bZ$-graded ring $A=\mathop{\oplus}\limits_{i\in \bZ} A_i$ is called $\bN$-graded if $A_i=0$ for all $i<0$. For a $\bZ$-graded ring $A=\mathop{\oplus}\limits_{i\in \bZ} A_i$, a left $A$-module $M$ is called graded if $M=\mathop{\oplus}\limits_{i\in \bZ} M_i$ for some subgroups $M_i$ such that $A_iM_j\subseteq M_{i+j}$ for all $i,j\in \bZ$.
The category of graded (left) $A$-modules is denoted by $\Gr A$, where $\Hom_{\Gr A}(M, N)=\{f: {}_AM \to {}_AN\mid f(M_i) \subseteq N_i,\, \forall \, i \in \bZ\}$, and $\gr A$ is the full subcategory of $\Gr A$ consisting of all finitely generated graded $A$-modules.

For a graded left $A$-module $M$ and $n\in \bZ$, let $M(n)$ be a graded left $A$-module with $i$-th degree part $M(n)_i=M_{n+i}$. For graded left $A$-modules $M$ and $N$, let
$$\gHom_A(M,N)=\mathop{\oplus}\limits_{n\in\bZ}\Hom_{\Gr A}(M,N(n)).$$
Let $\Ext_{\Gr A}^i(-,-)$ and $\gExt_A^i(-,-)$ be the $i$-th derived functors of $\Hom_{\Gr A}(-,-)$ and $\gHom_A(-,-)$ respectively.

For any $M\in \Gr A$, the graded projective dimension of $M$ is denoted by $\grpdim_AM$. The (left) graded global dimension of $A$ is denoted by $\grgldim A$.

The category of (left) $A$-modules is denoted by $\Mod A$. % and for convenience, we denote $\Hom_{\Mod A}(-,-)$ by $\Hom_A(-,-)$. Let $\Ext_A^i(-,-)$ be the $i$-th derived functor of $\Hom_A(-,-)$.
The projective dimension of an $A$-module $M$ is denoted by $\pdim_AM$. The (left) global dimension of $A$ is denoted by $\gldim A$.

The following lemma is \cite[Corollary A.I.2.7 and Theorem A.II.8.2]{NO1}
\begin{lemma}\label{grgldim and gldim}
    Let $A$ be a $\bZ$-graded ring.
    \begin{itemize}
        \item [(1)] $\grpdim_A M=\pdim_A M$ for any $M\in \Gr A$.
        \item [(2)] $\grgldim A\leqslant \gldim A\leqslant \grgldim A+1$.
    \end{itemize}
\end{lemma}

Let $M$ be a graded $A$-module. Then $M$ is noetherian in $\Gr A$ if and only if $M$ is noetherian in $\Mod A$ \cite[Theorem A.II.3.5]{NO1};  $M$ is projective in $\Gr A$ if and only if $M$ is projective in $\Mod A$ \cite[Corollary A.I.2.2]{NO1}. So, $A$ is graded (left or right) noetherian if and only if $A$ is (left or right) noetherian in ungraded sense.

\subsection*{Commutative graded rings}
In the rest of this section, we assume that $A$ is a commutative $\bZ$-graded ring.

A graded ideal $\fp$ of $A$ is called {\it graded prime} if it is a prime ideal of $A$. So, a graded ideal $\fp$ of $A$ is graded prime if and only if that $xy\in \fp$ implies $x\in \fp$ or $y\in \fp$ for any homogeneous elements $x,y\in A$.
Clearly, every maximal graded ideal is graded prime. The following lemma follows from \cite[Lemma 1.5.7]{BH}.

\begin{lemma}\label{structure of graded simple ring}
    Let $A$ be a $\bZ$-graded ring, $\fm$ be a maximal graded ideal of $A$. Let $k=(A/\fm)_0$. Then
    \begin{itemize}
        \item [(1)] $k$ is a field;
        \item [(2)] $A/\fm=k$ or $A/\fm=k[x,x^{-1}]$ for some homogeneous element $x$ of positive degree in $A/\fm$ which is transcendental over $k$.
    \end{itemize}
\end{lemma}
For any graded ideal $I$ of $A$, the \textit{graded height} $\grheight I$ of $I$ is defined as
$$ \max \{n \mid \exists \textrm{ a chain} \,
\fp_0\supsetneq \fp_1\supsetneq\cdots \supsetneq \fp_n \textrm{of graded prime ideals} \}$$
where $p_0$ runs over all minimal graded prime ideals containing $I$.

The \textit{graded Krull dimension} of $A$, denoted by $\grKdim A$, is defined to be the supremum of the graded heights.
For a $\bZ$-graded ring $A$, $\grKdim A$ may not equal to its Krull dimension $\Kdim A$. For example, if $A=k[x,x^{-1}]$ with the degree of $x$ being $1$, then $\grKdim A=0$ but $\Kdim A=1$.

If $A$ is noetherian and $\fp$ is graded prime, then $\grheight \fp=\height \fp$ (\cite[Theorem 1.5.8]{BH}). So for a graded prime ideal $\fp$, we will use $\height \fp$ to denote its graded height.

For any multiplicatively closed subset $S$ of $A$, let $S_h$ denote the subset of all homogeneous elements in $S$. Clearly, $S_h$ is also a multiplicatively closed subset. Obviously, if for any $s\in S$, there is at least one homogeneous component of $s$ contained in $S$, then $S^{-1}M = 0$
if and only if $S_h^{-1}M=0$ for any graded $A$-module $M$. Let
$$(S_h^{-1}M)_i=\{x/s\mid x\in M, s\in S \text{ are homogeneous such that}\, \deg x-\deg s=i\}.$$
Then $S_h^{-1}A=\mathop{\oplus}\limits_{i\in\bZ} (S_h^{-1}A)_i$ is a $\bZ$-graded ring and $S_h^{-1}M=\mathop{\oplus}\limits_{i\in\bZ} (S_h^{-1}M)_i$ is a graded $S_h^{-1}A$-module.
If $S=A\backslash \fp$ for a (graded) prime ideal $\fp$ of $A$,  then $M_{(\fp)}:=S_h^{-1}M$ is called the \textit{homogeneous localization} of $M$  at $\fp$.

Here are some properties of homogeneous localization which are used later.

\begin{lemma}\label{graded prime and ann and localization}
    Suppose $A$ is a $\bZ$-graded ring and $M$ is a finitely generated graded $A$-module.
    Then for any (graded) prime ideal $\fp$, $M_{(\fp)}\neq 0$ if and only if $\Ann_A(M)\subseteq \fp$; if and only if $M_{\fp}\neq 0$.
\end{lemma}

\begin{lemma}\label{localizaition preserves injective}
    Let $A$ be a $\bZ$-graded noetherian ring, and $E$ a graded injective $A$-module. Then, for any multiplicatively closed subset $S$ consisting of homogeneous elements of $A$, $S^{-1}E$ is a graded injective $S^{-1}A$-module.
\end{lemma}
\begin{proof}
    For any finitely generated graded $S^{-1}A$-module $L$, let $M$ be a finitely generated graded $A$-module such that $S^{-1}M=L$. Since $A$ is noetherian, $S^{-1}A$ is noetherian. By \cite[Corollary A.I.2.12]{NO1},
    \begin{align*}
        \gExt_{S^{-1}A}^1(L,S^{-1}E)&\cong \Ext_{S^{-1}A}^1(L,S^{-1}E)\\
        &\cong \Ext_A^1(M,E)\otimes_A S^{-1}A\\
        &\cong \gExt_A^1(M,E)\otimes_A S^{-1}A\\
        &=0.
    \end{align*}
   It follows from the graded version of Baer's theorem (\cite[Corollary 2.4.8]{NO2}) that $S^{-1}E$ is a graded injective $S^{-1}A$-module.
\end{proof}

\begin{definition}
    A $\mathbb{Z}$-graded ring $A$ is called \textit{graded local} if $A$ has only one maximal graded ideal.
\end{definition}

If $A$ is a ($\bZ$-graded) local ring, we usually use $\fm$ to denote the maximal ($\bZ$-graded) ideal of $A$ and write $k_A=A/\fm$. Sometimes, we will briefly say $(A,\fm,k_A)$ or $(A,\fm)$ is a ($\bZ$-graded) local ring.

\begin{example}
    (1) Let $\fp$ be a graded prime ideal of a $\bZ$-graded ring $A$. Then the homogeneous localization $A_{(\fp)}$ is a $\bZ$-graded local ring with maximal graded ideal $\fp A_{(\fp)}$.

    (2) If $A$ is an $\bN$-graded ring and $\fp$ is a graded prime ideal such that $A_{>0}$ is not contained in $\fp$, then $A_{(\fp)}$ is a $\bZ$-graded local ring with $(A_{(\fp)})_{<0}\neq 0$.
\end{example}

\section{Regular Graded Local Rings}\label{Graded Local Regular Rings}
In this section, we define and characterize regular graded local rings. All the rings considered in this section are  commutative.
Recall that  a noetherian local ring $(A,\fm,k_A)$ is \textit{regular} if $\dim_{k_A}(\fm/\fm^2)=\Kdim A$. 
There are lots of characterizations of regular noetherian local rings, for example, a famous result of Serre says that
$A$ is regular if and only if $\gldim A$ is finite, and in this case $\gldim A=\Kdim A$ (say, see \cite[Theorems 42 and 45]{Ma}). In general, a commutative noetherian ring  $A$ is called {\it regular} if $A_\mathfrak{p}$ is a regular local ring for any prime ideal $\mathfrak{p}$ of $A$.

\subsection{Characteristic polynomial of \texorpdfstring{$\bZ$}{Z}-graded local rings}
Suppose $(A,\fm,k_A)$ is a noetherian (resp. $\bZ$-graded) local ring. Let $l_A(M)$ (resp. $l^g_A(M)$) be the length of an $A$-module (resp. a graded $A$-module) $M$ of finite length.

\begin{lemma}\label{lengh and localization}
    Suppose $(A,\fm,k_A)$ is a $\bZ$-graded local ring. If $M$ is a graded $A$-module of finite length, then $M_{\fm}$ has finite length as an $A_{\fm}$-module, and $l^g_A(M)=l_{A_{\fm}}(M_{\fm})$.
\end{lemma}
\begin{proof}
    If $l^g_A(M)=1$, then $M\cong k_A(s)$ for some $s\in \bZ$. Since $(k_A)_{\fm}\cong A_{\fm}/\fm A_{\fm}$, $l^g_A(k_A)=l_{A_{\fm}}((k_A)_{\fm})=1$. So $l^g_A(M)=l_{A_{\fm}}(M_{\fm})=1$. An induction on the length of $M$ shows that $l^g_A(M)=l_{A_{\fm}}(M_{\fm})$.
\end{proof}

Suppose $(A,\fm,k_A)$ is a noetherian $\bZ$-graded local ring. A graded ideal $I$ satisfying that $\fm^s\subseteq I\subseteq \fm$ for some $s\in \bN$ is called a \textit{graded $\fm$-primary ideal}.
For any graded $\fm$-primary ideal $I$, $A/I^n$ is a graded $A$-module of finite length for any $n\in \bN$.

\begin{lemma}\label{chara poly}
    Suppose $(A,\fm,k_A)$ is a noetherian $\bZ$-graded local ring, and $I$ is a graded $\fm$-primary ideal generated by $m$ homogeneous elements.
    \begin{itemize}
        \item [(1)] There is a polynomial $\chi^g_I(t)\in \bQ[t]$ such that $l^g_A(A/I^n)=\chi^g_I(n)$ for $n\gg 0$.
        \item [(2)]  $\deg \chi^g_I(t) \leqslant m$.
    \end{itemize}
\end{lemma}
\begin{proof}
    (1) Since $I$ is a graded $\fm$-primary ideal, $I_{\fm}$ is an $\fm A_{\fm}$-primary ideal of the local ring $A_{\fm}$. Let $\chi_{I_{\fm}}(t)\in \bQ[t]$ be the characteristic polynomial of  $A_{\fm}$ relative to $I_{\fm}$.
    Then, $l_{A_{\fm}}(A_{\fm}/I^n_{\fm})=\chi_{I_{\fm}}(n)$ for $n \gg 0$. Since $A_{\fm}/I^n_{\fm}=(A/I^n)_{\fm}$, $l^g_A(A/I^n)=l_{A_{\fm}}(A_{\fm}/I^n_{\fm})$ by Lemma \ref{lengh and localization}. So $\chi_I^g(t)=\chi_{I_{\fm}}(t)$ is a polynomial we want to find.

    (2) By \cite[Proposition 1.5.15]{BH}, the minimal number of homogeneous generators of $I$ is equal to the minimal number of generators of $A_{\fm}$-module $I_{\fm}$. Since  $\deg \chi_{I_{\fm}}(t)$ is no more than the number of generators of $I_{\fm}$, $\deg \chi_I^g(t)\leqslant m$.
\end{proof}

The polynomial $\chi_I^g(t)$ is called the \textit{characteristic polynomial} of the graded local ring $A$ relative to $I$. By Lemma \ref{chara poly}, $\chi_I^g(t)=\chi_{I_{\fm}}(t)$, the characteristic polynomial of the local ring $A_{\fm}$ relative to $I_{\fm}$. If $Q$ is another graded $\fm$-primary ideal of $A$ , then $\chi_Q^g(t)=\chi_{Q_{\fm}}(t)$. Since $\deg \chi_{I_{\fm}}(t) = \deg \chi_{Q_{\fm}}(t)$, $\deg \chi^g_I(t) = \deg \chi^g_Q(t)$. The degree of $\chi^g_I(t)$ is independent of the choice of the graded $\fm$-primary ideals of $A$, which is denoted by $d(A)$.

\begin{proposition}
    Let $(A,\fm,k_A)$ be a noetherian $\bZ$-graded local ring. Then the following integers are equal.
    \begin{itemize}
        \item [(1)] $d(A)$. % the common degree of characteristic polynomials of graded ring $A$ relative to $\fm$-primary ideals.
        \item [(2)] $m(A)$, the minimal number of homogeneous elements generating a graded $\fm$-primary ideal.
        \item [(3)] $\grKdim A$, the graded Krull dimension of $A$.
    \end{itemize}
\end{proposition}
\begin{proof}
    Let $d(A_{\fm})$ be the degree of characteristic polynomial of $A_{\fm}$ relative to $\fm A_{\fm}$-primary ideals. Then $d(A_{\fm})=\Kdim A_{\fm}$. Since $\Kdim A_{\fm}=\height \fm=\grKdim A$, $d(A_{\fm})=\grKdim A$. It follows 
    that $d(A)=d(A_{\fm})=\grKdim A$.
    By Lemma \ref{chara poly}, $d(A)\leqslant m(A)$. To finish the proof, it suffices to prove that $m(A)\leqslant \grKdim A$.

    If $\grKdim A=0$, then $A$ is a graded artinian ring. So there is an integer $n$ such that $\fm^n=0$. Hence $0$ is a graded $\fm$-primary ideal. Then $m(A)=0$.

    Suppose $m(A)>0$. Let $\{P_1,\cdots,P_r\}$ be the set of all minimal prime ideals of $A$, all of which are graded prime. Clearly $\fm\nsubseteq P_i$ for all $i$. So $\fm \nsubseteq \cup P_i$. Take a homogeneous element $x\in \fm\backslash \cup P_i$. Then $(A/xA,\fm/x A)$ is a noetherian $\bZ$-graded local ring. Every chain of graded prime ideals in $A/xA$ is of the form
    $$P'_0/xA\supsetneq P'_1/xA\supsetneq \cdots \supsetneq P'_s/xA$$
    where $P'_i$ is a graded prime ideal of $A$ containing $xA$. Then there is some $1\leqslant i\leqslant r$ such that $P_i\subsetneq P'_s$. So $\grKdim A/xA+1\leqslant \grKdim A$. By induction hypothesis, $m(A/xA)\leqslant \grKdim A/xA$.

    On the other hand, every graded $(\fm/xA)$-primary ideal of $A/xA$ is of the form $Q/xA$ for some graded $\fm$-primary ideal $Q$ of $A$.
    Let $\{\bar{x}_1,\cdots,\bar{x}_{s}\}$ be a homogeneous generating subset of $Q/xA$ where $x_i\in Q$. Then $\{x_1,\cdots,x_s,x\}$ is a homogeneous generating subset of the $\fm$-primary ideal $Q$ of $A$. It follows that $m(A)\leqslant m(A/xA)+1$. So $m(A)\leqslant \grKdim A$.
\end{proof}

\subsection{Regular \texorpdfstring{$\bZ$}{Z}-graded local rings}
Since $k_A$ is a $\bZ$-graded simple ring for any  $\bZ$-graded local ring $(A,\fm,k_A)$, 
every graded $k_A$-module $M$ is a direct sum of graded simple $A$-modules, which are shifts of $k_A$ (\cite[Proposition 2.9.8]{NO2}). For a finitely generated graded $k_A$-module $M$, let $\rank_{k_A}(M)$ be the number of graded simple modules in its direct sum decomposition. Here is the definition of regular graded local ring.

\begin{definition}\label{def-gr-reg-local-ring}
    Let $(A,\fm,k_A)$ be a noetherian $\bZ$-graded local ring. If
    $$\rank_{k_A}(\fm/\fm^2)=\grKdim A,$$
    then $A$ is called a \textit{regular $\bZ$-graded local ring.}
\end{definition}

To characterize noetherian regular $\bZ$-graded local rings, we do some preparations.

    Let $A$ be a ring and $M$ a finitely generated $A$-module. Recall that an \textit{$M$-regular sequence} is a sequence $x_1,\cdots,x_n\in A$ such that $x_i$ is not a zero-divisor of $M/(x_1,\cdots,x_{i-1})M$ for $1\leqslant i\leqslant n$ and $M\neq (x_1,\cdots,x_n)M$.

The following lemma is \cite[Theorem 1.2.5]{BH}. 
\begin{lemma}\label{grade}
    Let $A$ be a noetherian ring, $I$ an ideal of $A$ and $M$ a finitely generated $A$-module such that $IM\neq M$.
    \begin{itemize}
        \item [(1)] All maximal $M$-regular sequences in $I$ have the same length. The common length of all maximal $M$-regular sequences in $I$ is denoted by $\grade(I,M)$, which is called the graded of $I$ on $M$.
        \item [(2)] $\grade(I,M)=\min\{i\mid \Ext_A^i(A/I,M)\neq 0\}.$
    \end{itemize}
\end{lemma}

Note that if $(A,\fm)$ is a $\bZ$-graded local ring and $M$ is a finitely generated graded $A$-module, then any homogeneous $M$-regular sequence is contained in $\fm$, because all the homogeneous elements in $A\backslash \fm$ are invertible.

Let $h(I)$ denote the set of all homogeneous elements of $I$ for any ideal $I$ of $A$.
\begin{definition} \label{completely-projective}\cite[B.III.3]{NO1}
    Let $A$ be a $\bZ$-graded ring. If, for any graded ideal $I$ and any finite set of graded prime ideals $P_1,\cdots,P_n$, $h(I) \subseteq P_1\cup\cdots \cup P_n$ implies that $I$ is contained in some $P_i$, then $A$ is called \textit{completely projective}.
\end{definition}

\begin{lemma}\label{graded-prime-avoidance}
    Let $(A, \fm, k_A)$ be a noetherian $\bZ$-graded local ring.  If $I$ a graded ideal of $A$ such that $h(I) \subseteq P_1\cup\cdots \cup P_n$ for graded prime ideals $P_1,\cdots,P_n$ not containing $A_{\geqslant 1} =\oplus_{i\geqslant 1} A_i$, then $I \subseteq P_i$ for some $1 \leqslant i \leqslant n$.
\end{lemma}
\begin{proof}
If some $P_i=\fm$, then $I \subseteq P_i$. Otherwise, the conclusion follows from the same proof as \cite[Lemma B.III.3.1]{NO1})
\end{proof}

It is easy to see that any $\bN$-graded ring $A$ such that $A_0$ is a field is completely projective.
Any $\bZ$-graded local ring $(A,\fm,k_A)$ such that $k_A$ is not a field is completely projective (see \cite[Example B.III.3.2]{NO1}). Next lemma is \cite[Corollary B.III.3.4]{NO1}.

\begin{lemma}\label{CP and lengh of homo regular sequence}
    Let $A$ be a noetherian $\bZ$-graded ring which is completely projective, $I$ a graded ideal of $A$ and $M$ a finitely generated graded $A$-module with $IM\neq M$. If $\grade(I,M)=n$, then there is an $M$-regular sequence in $I$ consisting of homogeneous elements with length being $\grade(I,M)$.
\end{lemma}

Typically, the assertion that any regular local ring of dimension $d$ possesses a regular sequence of length $d$ is proved via the prime avoidance lemma. In the context of a completely projective regular graded local ring of dimension $d$, the existence of such a sequence consisting of homogeneous elements can be proved by using Lemma \ref{CP and lengh of homo regular sequence}. Nevertheless, not all graded rings are completely projective, as evidenced by \cite[Example B.III.3.2]{NO1}. Consequently, an alternative approach is necessary to establish that any regular graded local ring of dimension $d$ contains a regular sequence of length $d$ consisting of homogeneous elements.

\begin{lemma}\label{rank and dim}
    Let $(A,\fm,k_A)$ be a noetherian $\bZ$-graded local ring and $\tilde{k}=A_{\fm}/\fm A_{\fm}$. Then
    $$\dim_{\tilde{k}}(\fm A_{\fm}/(\fm A_{\fm})^2)=\rank_{k_A}(\fm/\fm^2).$$
\end{lemma}
\begin{proof}
It follows from \cite[Proposition 1.5.15(a)]{BH} and Lemma \ref{lengh and localization}.
\end{proof}

The associated graded ring $\Gr_{\fm} A=A/\fm \oplus \fm/\fm^2\oplus\cdots=\bigoplus_{n\in \bN} \fm^n/\fm^{n+1}$ of $(A,\fm,k_A)$ with respect to the filtration $\cdots \subset \fm^{n+1} \subset \fm^n\subset \fm^{n-1}\subset \cdots$ is a bigraded ring, where the first grading is induced by the grading of $A$ and the second is induced by the filtration.

Now we are ready to characterize regular graded local rings.

\begin{theorem}\label{chara of graded regular}
    Let $(A,\fm,k_A)$ be a noetherian $\bZ$-graded local ring with graded Krull dimension $d$. Then the following statements 
    are equivalent.
    \begin{itemize}
        \item [(1)] $(A,\fm)$ is a regular graded local ring.
        \item [(2)] $(A_{\fm},\fm A_{\fm})$ is a regular local ring.
        \item [(3)] $\grgldim A$ is finite.
        \item [(4)] $\grgldim A=\grKdim A$.
        \item [(5)] $A$ is a regular ring.
        \item [(6)] $k_A[x_1,\cdots,x_d] \cong \Gr_{\fm} A$ as bigraded rings, where the bigrading of the polynomial ring $k_A[x_1,\cdots,x_d]$ is given as: the first grading of $k_A$ is inherited from $k_A=A/\fm$, and the second grading of elements in $k_A$ are zero; the second grading of the homogeneous elements $x_i$ are all $1$.
        \item [(7)] $A$ has a regular sequence of the length $\rank_{k_A}(\fm/\fm^2)$ consisting of homogeneous elements.
    \end{itemize}
\end{theorem}
\begin{proof}
    (1) $\Leftrightarrow$ (2) Let $\tilde{k}=A_{\fm}/\fm A_{\fm}$.
    Note that $\grKdim A=\height \fm=\Kdim A_{\fm}$. Then, by Lemma \ref{rank and dim},
    \begin{align*}
        &\rank_{k_A}(\fm/\fm^2)=\grKdim A\\
        \Leftrightarrow & \rank_{k_A}(\fm/\fm^2)=\Kdim A_{\fm}\\
        \Leftrightarrow & \dim_{\tilde{k}}(\fm A_{\fm}/(\fm A_{\fm})^2)=\Kdim A_{\fm}.
    \end{align*}
    So $(A,\fm)$ is a regular graded ring if and only if $(A_{\fm},\fm A_{\fm})$ is a regular local ring.

    (2) $\Rightarrow$ (3)  For any $M\in\gr A$, $\grpdim_AM=\pdim_{A_{\fm}}M_{\fm}$ by \cite[Proposition 1.5.15(e)]{BH}. So $\grpdim_AM \leqslant \gldim A_{\fm}$, which is finite as $A_{\fm}$ is a regular local ring.
    It follows that $\grgldim A$ is finite.

    (3) $\Rightarrow$ (4)
    Since $\grgldim A$ is finite, $\gldim A$ is finite by Lemma \ref{grgldim and gldim}. So, $\gldim A_{\fm}$ is finite, and consequently $A_{\fm}$ is a regular local ring. Hence $\gldim A_{\fm}=\Kdim A_{\fm}=\height \fm=\grKdim A$.

    Let $n=\grgldim A$. Then there is some $M\in\gr A$, such that $\grpdim_AM=n$. By \cite[Proposition 1.5.15(e)]{BH}, $\grpdim_AM=\pdim_{A_{\fm}}M_{\fm}=n$. Hence $\gldim A_{\fm}\geqslant n$. Therefore $\grKdim A\geqslant \grgldim A$.

    By \cite[Proposition 1.5.15(e)]{BH}, $\grade(\fm,A)=\grade(\fm A_{\fm},A_{\fm})$. By Lemma \ref{grade},
    $$\min\{i\mid \Ext_A^i(A/\fm,A)\neq 0\}=\min\{i\mid \Ext_{A_{\fm}}^i(A_{\fm}/\fm A_{\fm},A_{\fm})\neq 0\}.$$
    Since $A_{\fm}$ is a regular local ring, $A_{\fm}$ is a Cohen-Macaulay ring \cite[Corollary 2.2.6]{BH}. 
    Then $\Kdim A_{\fm}=\min\{i\mid \Ext_{A_{\fm}}^i(A_{\fm}/\fm A_{\fm},A_{\fm})\neq 0\}$.
    Note that $\Ext_A^i(A/\fm,A)=\gExt_A^i(A/\fm,A)$. Therefore,
    $$\grKdim A =\Kdim A_{\fm}=\min\{i\mid \gExt_A^i(A/\fm,A)\neq 0\}.$$
 It follows that $\grKdim A\leqslant \grgldim A$. Hence $\grKdim A=\grgldim A$.

   (4) $\Rightarrow$ (5) As $A$ is noetherian, $\grKdim A = \height \fm$ is finite. So,
    $\grgldim A$ is finite. Then by Lemma \ref{grgldim and gldim}, $\gldim A$ is finite.
  It follows that $A$ is a regular ring.

   (5) $\Rightarrow$ (2) It is direct from the definition.

   (1) $\Rightarrow$ (6) By definition, $\rank_{k_A}(\fm/\fm^2)=d$. Let $\{t_1,\cdots,t_d\}$ be a homogeneous generating set of $\fm$, such that $\fm/\fm^2=k_A\bar{t_1} + k_A\bar{t_2}+ \dots + k_A\bar{t_d}$. Consider the surjective morphism of the bigraded rings
   $$\varphi:k_A[x_1,\cdots,x_d]\to \Gr_{\fm}A,x_i\mapsto \bar{t_i}$$
   where $\bideg(x_i)=(\deg \bar{t_i}, 1)$, and $\bideg (\bar{a})=(\deg \bar{a}, 0)$ for any $\bar{a} \in k_A$.
   In fact, $\varphi$ is also a morphism of graded $A$-modules.

To prove that  $\varphi$ is an isomorphism, it suffices to prove that 
 $(k_A[x_1,\cdots,x_d])_{\fm}\cong (\Gr_{\fm}A)_{\fm}$.
 Note $(k_A[x_1,\cdots,x_d])_{\fm}\cong \tilde{k}[x_1,\cdots,x_d]$ and $(\Gr_{\fm}A)_{\fm}\cong \Gr_{\fm A_{\fm}}A_{\fm}$. 
 So there is a surjective morphism of graded rings:
   $$\varphi_{\fm}:\tilde{k}[x_1,\cdots,x_d]\to \Gr_{\fm A_{\fm}} A_{\fm},x_i \mapsto \overline{t_i/1}.$$
By (2), $(A_{\fm},\fm A_{\fm},\tilde{k})$ is a regular local ring of dimension $d$. So $\varphi_{\fm}$ is an isomorphism. It follows that $\varphi$ is an isomorphism.

   (6) $\Rightarrow$ (1) It follows from the isomorphism in (6) that $\rank_{k_A}(\fm/\fm^2)=d$.

   (6) $\Rightarrow$ (7) It follows from the isomorphism in (6) that $\rank_{k_A}(\fm/\fm^2)=d$.
   
   Let $\varphi:k_A[x_1,\cdots,x_d]\to \Gr_{\fm}A,x_i\mapsto \bar{t_i}$ be the graded isomorphism, where
   $t_1,\cdots,t_d$ are homogeneous elements of $\fm$ such that $\varphi(x_i)=\bar{t_i}\in \fm/\fm^2$.  
   
   If $t_1a=0$ for some homogeneous element $a\in \fm^s\backslash\fm^{s+1}$, then $\bar{t_1}\bar{a}=0\in \fm^{s+1}/\fm^{s+2} \subset \Gr_{\fm}A$ and thus $\bar{a}=0\in\fm^s/\fm^{s+1}$. It follows that $a\in\fm^{s+1}$, which is a contradiction. Hence $t_1 \in \fm$ is a regular element. 

   Let $\bar{A}=A/(t_1)$ and $\bar{\fm}=\fm/(t_1)$. Then $\bar{A}/\bar{\fm}\cong k_A$ and $(\bar{A},\bar{\fm},k_A)$ is a noetherian $\bZ$-graded local ring. We claim that $\Gr_{\bar{\fm}}\bar{A}\cong \Gr_{\fm}A/(\bar{t_1})$ as bigraded rings.
   
   With respect to the second degree, the $n$-th degree part of $\Gr_{\bar{\fm}}\bar{A}$ is isomorphic to $\fm^n/((t_1)\cap \fm^n+\fm^{n+1})$, and the $n$-th degree part of $\Gr_{\fm}A/(\bar{t_1})$ is isomorphic to $\fm^n/(t_1\fm^{n-1}+\fm^{n+1})$. To prove $\Gr_{\bar{\fm}}\bar{A}\cong \Gr_{\fm}A/(\bar{t_1})$, it suffices to show $t_1\fm^{n-1}= (t_1)\cap\fm^n$. Note that $t_1\fm^{n-1}\subseteq (t_1)\cap\fm^n$ is obvious. 

   Now suppose $t_1a\in (t_1)\cap\fm^n$, where $a\in \fm^l\backslash\fm^{l+1}$ is a homogeneous element. 
   Since $\Gr_{\fm} A$ is isomorphic to $k_A[x_1,\cdots,x_d]$, $\bar{t_1}$ is regular in $\Gr_{\fm}A$. Then $0 \neq \bar{t_1}\bar{a}\in \fm^{l+1}/\fm^{l+2}$. It follows that $t_1a\in \fm^{l+1}\backslash \fm^{l+2}$.  
   Since $t_1a\in \fm^n$, $n \leqslant l+1$. Hence $a\in \fm^l\subseteq \fm^{n-1}$, and $t_1a\in t_1\fm^{n-1}$. Therefore, $(t_1)\cap\fm^n \subseteq t_1\fm^{n-1}$.
   
   In conclusion $t_1\fm^{n-1}= (t_1)\cap\fm^n$ and $\Gr_{\bar{\fm}}\bar{A}\cong \Gr_{\fm}A/(\bar{t_1})$ as bigraded rings. 
   
   Then $\varphi$ induces an isomorphism of bigraded rings
   $$k_A[x_2,\cdots,x_d]\cong k_A[x_1,\cdots,x_d]/(x_1)\cong \Gr_{\fm}A/(\bar{t_1})\cong \Gr_{\bar{\fm}}\bar{A}.$$
   By a similar argument, $t_2+(t_1)\in\bar{A}$ is a regular element. So by induction, we have $\{t_1,\cdots,t_d\}$ is a regular sequence consisting of homogeneous elements of $A$. 

   (7) $\Rightarrow$ (2) 
   Let $n=\rank_{k_A}(\fm/\fm^2)$. Suppose $\{x_1,\cdots,x_n\}$ is a homogeneous $A$-regular sequence. Then $\{x_1/1,\cdots,x_n/1\}$ is an $A_{\fm}$-regular sequence \cite[Corollary 1.1.3]{BH}.
  Hence
   \begin{align*}
       n&\leqslant \grade(\fm A_{\fm},A_{\fm}) &\text{ (Lemma \ref{grade})}\\
        &\leqslant \Kdim A_{\fm} &\text{
 (\cite[Proposition 1.2.12]{BH})}\\
        &\leqslant \dim_{\tilde{k}}(\fm A_{\fm}/(\fm A_{\fm})^2) &\text{  (\cite[12.J]{Ma})}\\
        & =\rank_{k_A} (\fm/\fm^2)=n. &\text{  (Lemma \ref{rank and dim})}
   \end{align*}
   So, $\Kdim A_{\fm}=\dim_{\tilde{k}}(\fm A_{\fm}/(\fm A_{\fm})^2)$, that is, $(A_{\fm},\fm A_{\fm})$ is a regular local ring.
\end{proof}

\begin{corollary}
If $(A,\fm,k_A)$ is a noetherian regular $\bZ$-graded local ring, then $A$ is a domain.
\end{corollary}
\begin{proof} It follows from that $\bigcap_{n \in \bN} \fm^n =0$ and $\Gr_{\fm}A\cong k_A[x_1,\cdots,x_d]$.
\end{proof}

\section{Graded Isolated Singularities}\label{Graded Isolated Singularities}
In this section, we first assume that $A$ is a left noetherian (not necessarily commutative) $\bZ$-graded ring. We study graded isolated singularities from the perspective of projective schemes $\Proj A$ and $\qgr A$.

Let $J_A$ be the graded Jacobson radical of $A$, which is the intersection of all maximal graded left ideals of $A$. A graded $A$-module $M$ is called {\it torsion} if for any $x \in M$ there is $n \in \mathbb{N}$ such that $J_A^nx=0$.  If $M$ has no non-zero torsion submodule, then $M$ is called {\it torsion-free}.
Let $\Tor A$ (resp. $\tor A$) be the full subcategory of $\Gr A$ consisting of all (resp. finitely generated) torsion modules in $\Gr A$.
Note that $\tor A$ is a dense subcategory of $\gr A$, and $\Tor A$ is a localizing subcategory of $\Gr A$. Let
$$\QGr A=\Gr A/\Tor A \text{ and }\qgr A=\gr A/\tor A$$
be the quotient categories.
Then $\qgr A$ can be regarded as a full subcategory of $\QGr A$. For the theory of quotient categories, one can refer to \cite[Chapter 4]{Po}.

Let $\pi:\Gr A\to \QGr A$ be the quotient functor and $\omega$ be the right adjoint functor of $\pi$.
We write $\cM=\pi M$ for $M\in \Gr A$, and $\Hom_{\cA}(-,-)$  for the Hom functor in $\QGr A$.

Since $\QGr A$ is a Grothendieck category, it has enough injective objects. By \cite[Proposition 4.5.3]{Po}, the injective objects of $\QGr A$ are exactly the images of torsion-free graded injective $A$-modules in $\QGr A$. So, every object in $\QGr A$ has a minimal injective resolution. The length of the minimal injective resolution of $\cN\in \QGr A$ is called the {\it injective dimension} of $\cN$, denoted by $\idim_{\QGr A} \cN$.

The $i$-th right derived functor of $\Hom_{\cA}(\cM,-)$  is denoted  by $\Ext_{\cA}^i(\cM,-)$.
Then 
$$\idim_{\QGr A}(\cN)=\max\{i\mid \Ext_{\cA}^i(\cM,\cN)\neq 0 \textrm{ for some } \cM\in \QGr A\}.$$

\subsection{Ext groups in \texorpdfstring{$\qgr A$}{qgr A}}
Although $\qgr A$ may not have enough injective objects or projective objects, Ext groups in $\qgr A$ can be defined through its derived category. Let $D(\qgr A)$ be the derived category of $\qgr A$. Then the $i$-th derived functor of $\Hom_{\qgr A}(-,-)$ is defined by
$$\Ext_{\qgr A}^i(\cM,\cN):=\Hom_{D(\qgr A)}(\cM,\cN[i])$$
for $\cM,\cN\in \qgr A$, where $[i]$ is the $i$-th shift functor in $D(\qgr A)$.

In fact, for any $\cM,\cN\in \qgr A$ and $i\in \mathbb{N}$,
$$\Ext_{\qgr A}^i(\cM,\cN)\cong \Ext_{\cA}^i(\cM,\cN)$$
as showed in the next lemma. So, we may use the minimal injective resolution of $\cN$ in $\QGr A$ to compute the Ext group $\Ext_{\qgr A}^i(\cM,\cN)$.

As usual, $D^-(\qgr A)$ and $D^-(\QGr A)$ are the right bounded derived categories of $\qgr A$ and $\QGr A$ respectively, and $D^-_{\qgr A}(\QGr A)$ is the full subcategory of $D^-(\QGr A)$ consisting of the complexes whose cohomologies are in $\qgr A$.

\begin{lemma}\label{Ext in qgr and in QGr}
Let $A$ be a left noetherian $\bZ$-graded algebra. Then
$$D^-(\qgr A)\cong D^-_{\qgr A}(\QGr A).$$
In particular, for any $\cM,\cN\in \qgr A$ and $i\in\mathbb{N}$,
$$\Ext_{\qgr A}^i(\cM,\cN)\cong \Ext_{\cA}^i(\cM,\cN).$$
\end{lemma}
\begin{proof}
By a dual version of \cite[Proposition 1.7.11]{KS}, it suffices to prove that: if $\cM \to \cN$ is an epimorphism in $\QGr A$ with $\cN \in \qgr A$ then there is a morphism $\cL \to \cM$ with $\cL \in \qgr A$ such that the composition $\cL \to \cM \to \cN$ is epic.

Since $\cN \in \qgr A$, there exists a finitely generated graded module $N$ and a surjective morphism $M  \to N$ in $\Gr A$ such that $\pi M \cong \cM$ and $\pi N \cong \cN$ by \cite[Corollary 3.10]{Po}.
Then we may take a finitely generated graded submodule $L$ of $M$ so that
the restriction map $L \to N$ is surjective. Hence $\pi L = \cL \to \cN = \pi N$ is an epimorphism. It follows from the following diagram
$$
\xymatrix{
  L \ar@{^(->}[d] \ar[dr]^{}        \\
  M \ar[r]_{}  & N            }
$$
that the composition $\cL \to \cM \to \cN$ is an epimorphism.
\end{proof}

The global dimension $\gldim(\qgr A)$ is the global dimension of $\qgr A$ as an abelian category, that is,
  $$\gldim(\qgr A)=\max\{i\mid \Ext_{\qgr A}^i(\cM,\cN)\neq 0 \textrm{ for some } \cM,\cN\in \qgr A\}.$$

To study the Ext groups and the global dimension of $\qgr A$, the following lemmas are useful.

\begin{lemma}\label{idim cX}
Let $A$ be a left noetherian $\bZ$-graded ring. Suppose $\cQ, \cX \in\QGr A$.
\begin{itemize}
    \item [(1)] $\cQ$ is an injective object in $\QGr A$ if and only if $\Ext_{\cA}^1(\cM,\cQ)=0$ for any $\cM\in\qgr A$.
    \item [(2)] $\idim_{\QGr A}\cX=\max\{i\mid \Ext_{\cA}^i(\cM,\cX)\neq 0 \textrm{ for some } \cM\in \qgr A\}.$
    \item [(3)] $\gldim(\qgr A)=\max\{\idim_{\QGr A}\cN\mid\cN\in\qgr A\}.$
    \item [(4)] $\gldim(\qgr A)= \gldim(\QGr A)$.
\end{itemize}
\end{lemma}
\begin{proof}
    (1) One direction is clear. Suppose $\Ext_{\cA}^1(\cM,\cQ)=0$ for any $\cM\in \qgr A$. Let
    $0\to \cQ\to \cE^0\to \cE^1\to \cE^2\to \cdots$
    be the minimal injective resolution of $\cQ$. Then, for any $\cM\in \qgr A$, %the following sequence
    $$0\to \Hom_{\cA}(\cM,\cQ)\to \Hom_{\cA}(\cM,\cE^0)\to \Hom_{\cA}(\cM,\cE^1)\to \Hom_{\cA}(\cM,\cE^2)$$
    is exact. Hence, for any $M\in \gr A$,
    \begin{footnotesize}
    \begin{equation}\label{Hom(M,omega cE)}
        0\to \Hom_{\Gr A}(M,\omega\cQ)\to \Hom_{\Gr A}(M,\omega\cE^0)\to \Hom_{\Gr A}(M,\omega\cE^1)\to \Hom_{\Gr A}(M,\omega\cE^2)
    \end{equation}
    \end{footnotesize}
    is exact. By taking $M=A$, it follows that
    \begin{equation}\label{omega cE}
        0\to \omega \cQ\to \omega \cE^0\to \omega\cE^1 \to \omega\cE^2
    \end{equation}
    is exact. Since each $\omega \cE^i$ is a graded injective $A$-module \cite[Corollary 4.4.7]{Po}, the exact sequence (\ref{omega cE}) is a partial injective resolution of $\omega \cQ$. Hence, for any $M\in \gr A$, $\Ext_{\Gr A}^1(M,\omega\cQ)$ can be calculated by using the exact sequence (\ref{Hom(M,omega cE)}), and thus $\Ext_{\Gr A}^1(M,\omega\cQ)=0$. It follows that $\omega \cQ$ is injective. So $\cQ\cong \pi\omega \cQ$ is an injective object in $\QGr A$ by \cite[Proposition 4.5.3]{Po}.

    (2) It suffices to prove that if $\Ext_{\cA}^{n+1}(\cM,\cX)=0$ for any $\cM\in\qgr A$ then $\idim_{\QGr A}\cX\leqslant n$.
    Let
    $$0\to \cX\to \cE^0\to \cdots\to \cE^{n-1}\xrightarrow[]{d^n} \cE^n \to \cE^{n+1}\to \cdots$$
    be the minimal injective resolution of $\cX$.
    Let $\cK=\Coker d^n$. Then
    $$\Ext^1_{\cA}(\cM,\cK)\cong \Ext_{\cA}^{n+1}(\cM,\cX)=0$$
    for all $\cM\in\qgr A$. So $\cK$ is injective by (1). Thus $\idim_{\QGr A}\cX\leqslant n$.

    (3) It follows from (2) and Lemma \ref{Ext in qgr and in QGr}.

    (4) For any $\cN\in \QGr A$, $\omega\cN$ is a direct limit of finitely generated graded $A$-modules. Thus $\cN\cong\pi\omega \cN$ is a direct limit of objects $\cN_t$ in $\qgr A$. By a similar argument as \cite[Lemma 4.3.1]{BV}, for any $\cM\in \qgr A$ and $i\in \bN$, $\Ext_{\cA}^i(\cM,-)$ commutes with direct limits. Thus $\Ext_{\cA}^i(\cM,\cN)\cong \dlim\Ext_{\cA}^i(\cM,\cN_t)$. Then (4) follows from (2) and (3).
\end{proof}

\begin{lemma}\label{pi preserves injective hulls}
Let $A$ be a left noetherian $\bZ$-graded ring. If the largest torsion submodule of any graded injective $A$-module $I$ is a direct summand of $I$,
then the quotient functor $\pi:\Gr A\to \QGr A$ preserves injective hulls.
\end{lemma}
\begin{proof}
Let $X$ be a graded $A$-module and $I$ be the injective hull of $X$. Assume $I=I_1\oplus I_2$ where $I_1$ is torsion and $I_2$ is torsion-free. Since $\pi I_1=0$ and $\pi$ is exact, we have a monomorphism $\cX\to \pi I_2$ and $\pi I =\pi I_2$ is injective in $\QGr A$ by \cite[Proposition 4.5.3]{Po}. For any non-zero subobject $\cX'$ of $\pi I_2$, $\omega \cX'$ is a non-zero submodule of $\omega\pi I_2\cong I_2$.
 So, $\omega\cX'$ can be regarded as a submodule of $I_2$. Hence $\omega\cX'\cap X\neq 0$. Since $\pi$ is exact, $\pi(\omega\cX'\cap X)\subseteq \cX'\cap\cX$. Therefore $\cX'\cap\cX\neq 0$, and $\cX\to \pi I=\pi I_2$ is an essential extension. So, $\pi$ preserves injective hulls.
\end{proof}

\subsection{Graded isolated singularity}
 Now we are ready to define and characterize graded isolated singularities. First, we recall the definition of (graded) isolated singularities in commutative ($\bZ$-graded) local case.

\begin{definition}\label{def-iso-sing}
    Let $(A,\fm)$ be a commutative noetherian local ring. If $(A_{\fp},\fp A_{\fp})$ is regular for any non-maximal prime ideal $\fp$, then $A$ is called an \textit{isolated singularity}.

    If $(A,\fm)$ is a commutative noetherian $\bZ$-graded local ring and $(A_{(\fp)},\fp A_{(\fp)})$ is a regular graded local ring for any non-maximal graded prime ideal $\fp$, then $A$ is called a \textit{graded isolated singularity}.
\end{definition}

As it is well known that if $A$ is a graded quotient of polynomial rings then $\qgr A$ is equivalent to the category of the coherent sheaves over the projective scheme associated to $A$ \cite{Se}. This fact inspires the following definition by \cite{Jo,Ue1} in noncommutative projective geometry. 
We will show that $A$ is a graded isolated singularity if and only if the global dimension of $\qgr A$ is finite when $A$ is commutative.
Recall that (not necessarily commutative) $\bZ$-graded ring $A$ is called graded \textit{semilocal} if $A/J_A$ is a direct sum of left graded simple modules.

\begin{definition}\label{def-nc-iso-sing}
  Let $A$ be a left noetherian $\bZ$-graded semilocal ring.
  If $\gldim (\qgr A)$ is finite, then $A$ is called a \textit{noncommutative graded isolated singularity}.
\end{definition}

\begin{lemma}\label{max graded ideal and semilocal}
    If $A$ is a commutative $\bZ$-graded ring, then $A$ is graded semilocal if and only if $A$ has only finitely many maximal graded ideals.
\end{lemma}

\begin{corollary}\label{grKdim of comm semilocal ring}
  Any commutative noetherian $\bZ$-graded semilocal ring has finite graded Krull dimension.
\end{corollary}

\begin{lemma}\label{torsion and localization}
    Let $A$ be a commutative noetherian $\bZ$-graded semilocal ring and $M$ a graded $A$-module. Then $M$ is torsion if and only if $M_{(\fp)}=0$ for any non-maximal graded prime ideal $\fp$.
\end{lemma}
\begin{proof}
    Suppose $M$ is torsion. If there is a graded prime ideal $\fp$ such that $M_{(\fp)}\neq 0$, then there is some homogeneous element $x \in M$ such that $0 \neq x/1 \in M_{(\fp)}$. It follows that
     $J_A^r \subseteq \Ann_A(x) \subseteq \fp$ for some $r$.
    Hence $J_A\subseteq \fp$. Since $A$ is graded semilocal, by Lemma \ref{max graded ideal and semilocal}, $A$ has only finitely many maximal graded ideals, say, $\fm_1,\cdots,\fm_s$. Then
    $$\fm_1\cdots \fm_s\subseteq\fm_1\cap \cdots \cap\fm_s = J_A \subseteq \fp.$$
    Hence $\fm_i\subseteq \fp$ for some $i$, and $\fp=\fm_i$ is maximal.

    Conversely, suppose $M_{(\fp)}=0$ for any non-maximal graded prime ideal $\fp$. Then, for any homogeneous element $x\in M$,
    $\Ann_A(x)\nsubseteq \fp$.
    Let $\fm_1,\cdots, \fm_t,\fm_{t+1},\cdots,\fm_s$ be the set of all maximal graded ideals of $A$ such that $\Ann_A(x)$ is contained in $\fm_i$ only when $1\leqslant i\leqslant t$.  So, $\fm_1/\Ann_A(x),\cdots,\fm_t/\Ann_A(x)$, which are maximal, are exactly all the graded prime ideals of $A/\Ann_A(x)$. By \cite[Proposition 2.11.1]{NO2},
    the intersection of all the maximal graded ideals of $A/\Ann_A(x)$ is contained in the intersection of all the prime ideals of $A/\Ann_A(x)$, which is nilpotent. So $(\fm_1\cap \cdots\cap \fm_t)^d\subseteq \Ann_A(x)$ for some $d\in \bN$. It follows from $J_A\subseteq \fm_1\cap \cdots \cap \fm_t$ that $(J_A)^d\subseteq \Ann_A(x)$. Thus $(J_A)^d x=0$. So, $M$ is torsion.
\end{proof}

\begin{corollary}\label{torsion and injective hull}
    Let $A$ be a commutative noetherian $\bZ$-graded semilocal ring. Let $\fp$ be a graded prime ideal of $A$, and $E(A/\fp)$ be the graded injective hull of ${}_A(A/\fp)$.
    \begin{itemize}
        \item [(1)] $E(A/\fp)$ is torsion if and only if $\fp$ is a maximal graded ideal.
        \item [(2)] $E(A/\fp)$ is torsion-free if and only if $\fp$ is not a maximal graded ideal.
    \end{itemize}
\end{corollary}
\begin{proof} 
    (1) For any graded prime ideal $\mathfrak{q}$, $E(A/\fp)_{\mathfrak{q}}=0$ if and only if $\fp\nsubseteq \mathfrak{q}$. Hence $E(A/\fp)_{(\mathfrak{q})}=0$ if and only if $\fp\nsubseteq \mathfrak{q}$. It follows from Lemma \ref{torsion and localization} that $E(A/\fp)$ is torsion if and only if $\fp$ is a maximal graded ideal.

   (2) Suppose $\fp$ is a graded prime ideal but not maximal. If there is an element $0 \neq x\in A/\fp$ such that $(J_A)^nx=0$ for some $n\in \bN$, then $(J_A)^n\subseteq \Ann_A(x)\subseteq \fp$. So, $J_A\subseteq \fp$. Let $\fm_1,\cdots,\fm_s$ be the set of all maximal graded ideals of $A$. Then $\fm_1\cdots\fm_s\subseteq \fp$, which is a contradiction. Therefore $A/\fp$ is torsion-free. Since $E(A/\fp)$ is an essential extension of $A/\fp$, $E(A/\fp)$ is torsion-free.

   If $E(A/\fp)$ is torsion-free, then $\fp$ is not a maximal graded ideal by (1).
\end{proof}

\begin{corollary} \label{pi preserves minimal}
    If $A$ is a commutative noetherian $\bZ$-graded semilocal ring, then any graded injective $A$-module $I$ can be decomposed into $I_1\oplus I_2$ where $I_1$ is torsion and $I_2$ is torsion-free.

    Consequently, $\pi:\Gr A\to \QGr A$ preserves injective hulls.
\end{corollary}
\begin{proof} Any graded injective $A$-module
    $I$ can be decomposed into a direct sum of graded indecomposable injective modules, and every graded indecomposable injective module is of the form $E(A/\fp)(n)$ where $\fp$ is a graded prime ideal, $E(A/\fp)$ is the graded injective hull of $A/\fp$ and $n$ is an integer (see \cite[Theorem 3.6.3(b,c)]{BH}).

    By Corollary \ref{torsion and injective hull}, each indecomposable direct summands of $I$ is either torsion or torsion-free. Let $I_1$ be the direct sum of indecomposable torsion direct summands and $I_2$ be the direct sum of indecomposable torsion-free direct summands in the decomposition of $I$. Then $I=I_1\oplus I_2$ is the desired decomposition.
    The last statement follows from Lemma \ref{pi preserves injective hulls}.
\end{proof}

\begin{theorem}\label{char-graded-isolated-singualrty}
    Let $A$ be a commutative noetherian $\bZ$-graded semilocal ring with $\grKdim A=d$.
    Then the following are equivalent.
    \begin{itemize}
        \item [(1)] $(A_{(\fm)},\fm A_{(\fm)})$ is a graded isolated singularity for any maximal graded ideal $\fm$ of $A$.
        \item [(2)] $(A_{\fm},\fm A_{\fm})$ is an isolated singularity for any maximal graded ideal $\fm$ of $A$.
        \item [(3)] $(A_{(\fp)},\fp A_{(\fp)})$ is a regular graded local ring for any non-maximal graded prime ideal $\fp$ of $A$.
        \item [(4)] $(A_{\fp},\fp A_{\fp})$ is a regular local ring for any non-maximal graded prime ideal $\fp$ of $A$.
        \item [(5)] The global dimension of $\qgr A$ is $d-1$.
        \item [(6)] The global dimension of $\qgr A$ is finite.
    \end{itemize}
\end{theorem}
\begin{proof}
    Let $M$ be a finitely generated graded $A$-module, and $0\to M\to I^0\to \cdots\to I^{d-1}\to I^d\to \cdots$ be the minimal graded injective resolution of $M$.

    By Corollary \ref{pi preserves minimal},
$$0\to \pi M\to \pi I^0\to \cdots\to \pi I^{d-1}\to \pi I^d\to \cdots$$
    is the minimal injective resolution of $\pi M$ in $\QGr A$.

By Lemma \ref{localizaition preserves injective}, for any non-maximal prime ideal $\mathfrak{p}$,
$$0\to M_{(\fp)}\to I_{(\fp)}^0\to \cdots\to I_{(\fp)}^{d-1}\to I_{(\fp)}^d\to \cdots$$
is a graded injective resolution of $M_{(\fp)}$, which is minimal by \cite[Lemma A.I.2.8]{NO1} and \cite[Corollary 1.3]{Ba}.

    (1) $\Rightarrow$ (2)  For any prime ideal $\fp$ of $A$ properly contained in $\fm$, let $\fp^*$ be the graded ideal generated by homogeneous elements of $\fp$, which is a graded prime ideal of $A$, see for example \cite[Lemma 1.5.6]{BH}.
    Since $A_{(\fp^*)}$ is the homogeneous localization of $A_{(\fm)}$ at the non-maximal graded prime ideal
    $\fp^*A_{(\fm)}$, $(A_{(\fp^*)},\fp^* A_{(\fp^*)})$ is a regular graded local ring. It follows from Theorem \ref{chara of graded regular} and \cite[Exercise 2.2.24]{BH} that $(A_{\fp},\fp A_{\fp})$ is regular. Therefore $(A_{\fm},\fm A_{\fm})$ is an isolated singularity.

    (2) $\Rightarrow$ (1) For any graded prime ideal $\fp$ of $A$ properly contained in $\fm$, $(A_{\fp},\fp A_{\fp})$ is a regular local ring. Note that $(A_{\fp},\fp A_{\fp})$ is the (non-homogeneous) localization of $A_{(\fp)}$ at $\fp A_{(\fp)}$. By (2) $\Rightarrow$ (1) in Theorem \ref{chara of graded regular}, $(A_{(\fp)}, \fp A_{(\fp)})$ is a regular graded local ring. So $(A_{(\fm)},\fm A_{(\fm)})$ is a graded isolated singularity.

    (1) $\Leftrightarrow$ (3) By the definition of graded isolated singularity.

    (3) $\Leftrightarrow$ (4) It follows from Theorem \ref{chara of graded regular}.

    (3) $\Rightarrow$ (5)  For any non-maximal graded prime ideal $\mathfrak{p}$, the height of $\mathfrak{p}$ is less than $d$. So $\gldim A_{(\mathfrak{p})}<d$ by Theorem \ref{chara of graded regular}. It follows that $I^d_{(\fp)}=0$. By Lemma \ref{torsion and localization}, $I^d$ is torsion, so $\pi I^d=0$. Hence the injective dimension of $\pi M$ is no more than $d$. Therefore, $\gldim \qgr A<d$  by Lemma \ref{idim cX}.

Let $\mathfrak{q}$ be a graded prime ideal of height $d-1$. Then $\grgldim A_{(\mathfrak{q})}=d-1$ by Theorem \ref{chara of graded regular}. So there is a finitely generated graded $A$-module $N$ such that $(I_N^{d-1})_\mathfrak{q}\neq 0$, where $I_N^{d-1}$ is the $(d-1)$-th term in the minimal graded injective resolution of $N$. It follows that $\pi I^{d-1}_N\neq 0$. By Lemma \ref{idim cX}, $\gldim \qgr A=d-1$.

(5) $\Rightarrow$ (6) Obviously.

(6) $\Rightarrow$ (3) Suppose $\gldim \qgr A=l$ is finite. For any finitely generated graded $A_{(\fp)}$-module $L$, there is a finitely generated graded $A$-module $M$ such that $M_{(\fp)}=L$. Take the minimal graded injective resolution of $M$ as in the beginning of the proof. Then, for $i>l$, $\pi I^i=0$ by Lemma \ref{idim cX}, and so $I^i_{(\fp)}=0$ by Lemma \ref{torsion and localization}. Hence the graded injective dimension of $L$ is no more than $l$. So $\grgldim A_{(\fp)} \leqslant l$. By Theorem \ref{chara of graded regular}, $(A_{(\fp)},\fp A_{(\fp)})$ is a regular graded local ring.
\end{proof}

Next corollary justifies the definition of noncommutative isolated singularities (see Definition \ref{def-nc-iso-sing}).

\begin{corollary}\label{chara of isolated singularity}
    Let $(A,\fm,k_A)$ be a commutative noetherian $\bZ$-graded local ring with $\grKdim A=d$. Then the following are equivalent.
    \begin{itemize}
        \item [(1)] $(A,\fm)$ is a graded isolated singularity.
        \item [(2)] $(A_{\fm},\fm A_{\fm})$ is an isolated singularity.
        \item [(3)] The global dimension of $\qgr A$ is $d-1$.
        \item [(4)] The global dimension of $\qgr A$ is finite.
    \end{itemize}
\end{corollary}

    Let $k[x_1,\cdots,x_n]$ be the graded polynomial algebra over a field $k$ with $\deg x_i=1$. Let $A=k[x_1,\cdots,x_n]/I$ where $I$ is a proper graded ideal of $k[x_1,\cdots,x_n]$. Then $A$ is a commutative noetherian $\bN$-graded local ring with maximal graded ideal $\fm=A_{>0}$. In particular, for any non-maximal graded prime ideal $\fp$ of $A$, $(A_{(\fp)},\fp A_{(\fp)})$ is completely projective (\cite[Example B.III.3.2]{NO1}).
    Let $X=\Proj A$ be the projective scheme associated to $A$ and $\coh X$ be the category of coherent sheaves over $X$. Then $\coh(X)\cong \qgr A$ (\cite{Se}).

    As a corollary, we show that $(A,\fm)$ is a graded isolated singularity if and only if $\Proj A$ is smooth. Recall that a scheme $X$ is said to be \textit{smooth (or nonsingular)} if for every point $x\in X$, the stalk $\mathcal{O}_x$ is a regular local ring \cite{Ha}.

\begin{corollary}\label{iso-sing-of-quo-of-poly}
Let $A=k[x_1,\cdots,x_n]/I$ be a graded quotient of the polynomial algebra with $\deg x_i=1$ and $\fm=A_{>0}$. Let $\Proj A$ be the projective scheme associated to $A$. Then the following are equivalent.
\begin{itemize}
    \item [(1)] $(A,\fm)$ is a graded isolated singularity.
    \item [(2)] $(A_{\fm},\fm A_{\fm})$ is an isolated singularity.
    \item [(3)] The global dimension of $\qgr A$ is finite.
    \item [(4)] The global dimension of $\coh (\Proj A)$ is finite.
    \item [(5)] For any $\fp\in\Spec A\backslash\{\fm\}$, $(A_{\fp},\fp A_{\fp})$ is a regular local ring.
    \item [(6)] $\Proj A$ is smooth.
\end{itemize}
In this case, the global dimensions of $\qgr A$ and $\coh (\Proj A)$ are $\grKdim A-1$.
\end{corollary}
\begin{proof}
    The equivalences of (1)-(4) follow from Corollary \ref{chara of isolated singularity} and $\coh(X)\cong \qgr A$.

    (1) $\Rightarrow$ (5) For any $\fp$ not equal to $\fm$, $\fp^*$ is the graded ideal generated by the homogeneous elements of $\fp$, which is a prime ideal. Let $S=A\backslash\fp$ and $\tilde{S}=A\backslash\fp^*$. Then $S_h=\tilde{S}_h$. Since $\tilde{S}_h\subseteq S$, $A_{\fp}$ is a localization of $A_{(\fp^*)}$. Hence
    $\gldim A_{\fp}\leqslant \gldim A_{(\fp^*)}\leqslant \grgldim A_{(\fp^*)}+1$ by Lemma \ref{grgldim and gldim}. Since $A_{(\fp^*)}$ is regular, $\grgldim A_{(\fp^*)}$ is finite. Therefore $(A_{\fp},\fp A_{\fp})$ is a regular local ring. 

    (5) $\Rightarrow$ (2) is trivial.
    
    (1) $\Leftrightarrow$ (6) By \cite[Proposition II.2.5]{Ha}, for every $\fp\in \Proj A$, the stalk $\mathcal{O}_{\fp}$ is isomorphic to $(A_{(\fp)})_0$. So $\Proj A$ is smooth if and only if for any non-maximal graded prime ideal $\fp$ of $A$, $(A_{(\fp)})_0$ is a regular local ring.

   Let $B=A_{(\fp)}$ and $S_h$ be the set of homogeneous elements in $A\backslash \fp $ for a fixed non-maximal graded prime ideal $\fp$. We claim that $B$ is strongly graded, that is, for any $i,j\in \bZ$, $B_iB_j=B_{i+j}$. Since $\fp$ is not maximal, there is some $f\in A_1\backslash \fp$. Then for any integer $n$ and for any homogeneous element $g\in B_n$, $g=f^n(f^{-n}g)$. It follows that $B=B_0[f,f^{-1}]$. So $B$ is strongly graded.

    Since $B$ is strongly graded, the categories $\Gr B$ and $\Mod B_0$ are equivalent by \cite[Theorem A.I.3.4]{NO1}, and thus
    $\grgldim B=\gldim B_0$. It follows from Theorem \ref{chara of graded regular} and \cite[Theorem 42]{Ma} that $B$ is graded regular if and only if $B_0$ is regular. So $\Proj A$ is smooth if and only if $(A,\fm)$ is a graded isolated singularity.
\end{proof}

\section*{Acknowledgements} The authors are very grateful to the referee for his very valuable comments and suggestions, and for pointing out a mistake in an example in the previous version.

\end{document}